\documentstyle[amsfonts, amssymb, latexsym,12pt]{amsart}

\newtheorem{theorem}{Theorem}[section]
\newtheorem{lemma}[theorem]{Lemma}
\newtheorem{proposition}[theorem]{Proposition}
\newtheorem{corollary}[theorem]{Corollary}
\newtheorem{definition}[theorem]{Definition}

\newtheorem{question}[theorem]{Question}

\def\C{{\mbox{\rm\kern.24em
\vrule width.03em height1.43ex depth-.052ex \kern-.26em C}}}
\def\QSet{\mbox{\rm\kern.24em
\vrule width.03em height1.48ex depth-.051ex \kern-.26em Q}}
\def\Z{{\bf Z}}
\def\R{{\mbox{\rm I\kern-.22em R}}}
\def\P{{\bf P}}

\def\D{{\bf D}}
\def\T{{\bf T}}
\def\supp{{\rm supp}}
\def\size{{\rm size}}

\def\energy{{\rm energy}}
\def\sgn{{\rm sgn}}

\def\L{{\cal{L}}}
\def\D{{\cal{D}}}
\def\S{{\cal{S}}}

\def\dist{{\rm dist}}

\def\111{\gamma}

\def\be#1{\begin{equation}\label{#1}}
\def\bas{\begin{align*}}
\def\eas{\end{align*}}
\def\bi{\begin{itemize}}
\def\ei{\end{itemize}}
\newenvironment{proof}{\noindent {\bf Proof} }{\endprf\par}
\def \endprf{\hfill  {\vrule height6pt width6pt depth0pt}\medskip}
\def\emph#1{{\it #1}}

\title{Bi-parameter paraproducts}

\author{Camil Muscalu}
\address{Department of Mathematics, Cornell University, Ithaca, NY 14853}
\email{camil@@math.cornell.edu}

\address{Current Address: School of Mathematics, IAS,  Princeton, NJ 08540}
\email{camil@@math.ias.edu}

\author{Jill Pipher}
\address{Department of Mathematics, Brown University, Providence, RI 02912}
\email{jpipher@@math.brown.edu}

\author{Terence Tao}
\address{Department of Mathematics, UCLA, Los Angeles, CA 90095}
\email{tao@@math.ucla.edu}

\author{Christoph Thiele}
\address{Department of Mathematics, UCLA, Los Angeles, CA 90095}
\email{thiele@@math.ucla.edu}

\begin{document}

\begin{abstract}
In the first part of the paper we prove a bi-parameter version of a well known multilinear theorem of Coifman and Meyer.
As an example of an application of our main theorem, we generalize the Kato-Ponce inequality in nonlinear PDE. Then, we show that the double bilinear Hilbert transform does not satisfy any $L^p$ estimates.
\end{abstract}

\maketitle

\section{Introduction}

Let $f\in\cal{S}$$(\R^{2})$ be a Schwartz function in the plane. A well known inequality in elliptic PDE says that
\begin{equation}\label{ell}
\|\frac{\partial^2 f}{\partial x_1 \partial x_2}\|_p \lesssim
\|\Delta f\|_p
\end{equation}
for $1<p<\infty$, where $\Delta=\frac{\partial^2}{\partial x_1^2} + \frac{\partial^2}{\partial x_2^2}$ is the Laplace operator.

To prove (\ref{ell}) one just has to observe that
$$\frac{\partial^2 f}{\partial x_1 \partial x_2} = c R_1 R_2\Delta f$$
where

$$R_j f(x) = \int_{\R^2}\frac{\xi_j}{|\xi|} \widehat{f}(\xi) e^{2\pi i x\xi}\,d\xi$$
$j=1,2$ are the Riesz transforms and they are bounded linear operators on $L^p(\R^2)$ \cite{stein}.

An estimate of a similar flavour in non-linear PDE is the following inequality of Kato and Ponce \cite{kp}.
If $f,g\in\cal{S}$$(\R^2)$ and $\widehat{\D^{\alpha}f}(\xi):= |\xi|^{\alpha}\widehat{f}(\xi)$ $\alpha >0$, is the homogeneous
derivative, then
\begin{equation}\label{katoponce}
\|\D^{\alpha}(fg)\|_r \lesssim \|\D^{\alpha}f\|_p \|g\|_q + \|f\|_p \|\D^{\alpha}g\|_q
\end{equation}
for $1<p,q, r < \infty$ with the property $1/r=1/p+1/q$.

Heuristically, if $f$ oscillates more rapidly than $g$, then $g$ is essentially constant with respect to $f$ and so
$\D^{\alpha}(fg)$ behaves like $(\D^{\alpha}f)g$. Similarly, if $g$ oscillates more rapidly then $f$ then one expects
$\D^{\alpha}(fg)$ to be like $f(\D^{\alpha}g)$ and this is why there are two terms on the right hand side of
(\ref{katoponce}). In order to make this argument rigorous, one needs to recall the classical Coifman-Meyer theorem \cite{cm}, \cite{gt}, \cite{ks}.
Let $m$ be a bounded function on $\R^4$, smooth away from the origin and satisfying
\begin{equation}\label{1par}
|\partial^{\beta}m(\gamma)| \lesssim \frac{1}{|\gamma |^{|\beta|}}
\end{equation}
for sufficiently many $\beta$. Denote by $T_m(f, g)$ the bilinear operator defined by
\begin{equation}\label{defop}
T_m(f, g)(x) = \int_{\R^4} m(\xi, \eta) \widehat{f}(\xi) \widehat{g}(\eta) e^{2\pi i x(\xi +\eta)} \, d\xi d\eta.
\end{equation}
Then, $T_m$ maps $L^p\times L^q\rightarrow L^r$ as long as $1<p,q\leq \infty$, $1/r=1/p+1/q$ and $0<r<\infty$.

This operator takes care of the inequality (\ref{katoponce}) in essentially the same way in which the Riesz transforms
take care of (\ref{ell}). The details will be presented later on in the Appendix (see also \cite{kp}).

But sometimes (see \cite{kpv}), in non-linear PDE one faces the situation when a partial differential operator such as
$\widehat{\D_1^{\alpha}\D_2^{\beta}f}(\xi_1, \xi_2):=|\xi_1|^{\alpha} |\xi_2|^{\beta} \widehat{f}(\xi_1, \xi_2)$ $\alpha$, $\beta >0$, acts on a nonlinear
expression such as the product of two functions. It is therefore natural to ask if there is an inequality analogous
to (\ref{katoponce}) for these operators. The obvious candidate, according to the same heuristics, is the following
inequality.

\begin{equation}\label{2katoponce}
\|\D_1^{\alpha}\D_2^{\beta}(fg)\|_r \lesssim  \|\D_1^{\alpha}\D_2^{\beta}f\|_p \|g\|_q + \|f\|_p \|\D_1^{\alpha}\D_2^{\beta}g\|_q + 
 \|\D_1^{\alpha}f\|_p \|\D_2^{\beta}g\|_q + \|\D_1^{\alpha}g\|_p \|\D_2^{\beta}f\|_q.
\end{equation}   
If one tries to prove it, one realizes that one needs to understand bilinear operators whose symbols satisfy estimates of the form
\begin{equation}\label{2par}
\left|\partial^{\alpha_1}_{\xi_1}
\partial^{\alpha_2}_{\xi_2}
\partial^{\beta_1}_{\eta_1}
\partial^{\beta_2}_{\eta_2} m(\xi, \eta)\right|\lesssim
\frac{1}{|(\xi_1, \eta_1)|^{\alpha_1 +\beta_1}}
\frac{1}{|(\xi_2, \eta_2)|^{\alpha_2 +\beta_2}}.
\end{equation}
Clearly, the class of symbols verifying (\ref{2par}) is strictly wider then the class of symbols satisfying (\ref{1par}).
These new $m$'s behave as if they were products of two homogeneous symbols of type (\ref{1par}), one of variables
$(\xi_1,\eta_1)$ and the other of variables $(\xi_2,\eta_2)$.

The main task of the present paper is to prove $L^p$ estimates for such operators in this more delicate product setting.
Our main theorem is the following.

\begin{theorem}\label{main}
If $m$ is a symbol in $\R^4$ satisfying (\ref{2par}), then the bilinear operator $T_m$ defined by (\ref{defop})
maps $L^p\times L^q\rightarrow L^r$ as long as $1<p,q\leq \infty$, $1/r=1/p+1/q$ and $0<r<\infty$.
\end{theorem}

It will be clear from the proof of the theorem that the $n$-linear analogue of this result is also true (see Section 8 for a precise statement).
Particular cases of this theorem have been considered by Journ\'e (see \cite{journe} and also \cite{cj}) who proved that in the situation of tensor products of two generic paraproducts, one has
$L^2\times L^{\infty}\rightarrow L^2$ estimates. Our approach is different from his and is based on arguments
with a strong geometric structure. The reader will notice that part of the difficulties of the general case comes from the fact that there is no analogue of the classical Calder\'on-Zygmund
decomposition in this bi-parameter framework and so the standard argument \cite{cm}, \cite{ks}, \cite{gt} used to prove such estimates, has to be changed.

The paper is organized as follows. In the next section, we discretize our operator and reduce it to a biparameter general
paraproduct. In the third section we present a new proof of the classical one parameter case. This technique will be very helpful to handle an error term later on 
in section six. Sections four, five and six are devoted to the proof of our main theorem (\ref{main}).
 Section seven contains a counterexample to the boundedness of the double bilinear Hilbert transform and then, the paper ends with some further comments and open questions. In the Appendix we explain how theorem \ref{main} implies
 inequality (\ref{2katoponce}).

{\bf Acknowledgements}: We would like to express our thanks to Carlos Kenig for valuable conversations and to the referees for their suggestions,
 which improved the presentation of the paper.

The first two authors were partially supported by NSF Grants. The third author is a Clay Prize Fellow and is partially supported by a
Packard Foundation Grant. The fourth author was partially supported by the NSF Grants DMS 9985572 and DMS 9970469.

\section{Reduction to bi-parameter paraproducts}

In order to understand the operator $T_m$, the plan is to carve it into smaller pieces well adapted to its biparameter structure. First, by writing
the characteristic functions of the planes $(\xi_1, \eta_1)$ and $(\xi_2, \eta_2)$ as finite sums of smoothed versions of characteristic functions of cones of the form
$\{(\xi,\eta) : |\xi|\leq C|\eta|\}$ or $\{(\xi,\eta) : |\xi|\geq C|\eta|\}$  , we decompose our operator into a finite sum of several parts. 
Since all the operators obtained in this decomposition can be treated in the same way, we will discuss in 
detail only one of them, which will be carefully defined below (in fact, as the reader will notice, the only difference between any arbitrary case and the one we will explain
here, is that the functions $MM$, $SS$, $MS$, $SM$ defined later on at page 6, have to be moved around).

Let $\phi$, $\psi$ be two Schwartz bumps on $[0,1]$, symmetric with respect to the origin and such that
$\supp(\widehat{\phi})\subseteq [-1/4,1/4]$ and $\supp(\widehat{\psi})\subseteq [3/4,5/4]$. Recall the translation
and dilation operators $\tau_h$, $D^p_{\lambda}$ given by
$$\tau_hf(x)=f(x-h)$$
$$D^p_{\lambda}f(x)=\lambda^{-1/p}f(\lambda^{-1}x)$$
and then define

$$C'(\xi_1, \eta_1)=\int_{\R} D_{2^{k'}}^{\infty}\widehat{\phi}(\xi_1) D_{2^{k'}}^{\infty}\widehat{\psi}(\eta_1)\,dk'$$
and

$$C''(\xi_2, \eta_2)=\int_{\R} D_{2^{k''}}^{\infty}\widehat{\psi}(\xi_2) D_{2^{k''}}^{\infty}\widehat{\phi}(\eta_2)\,dk''.$$
As we said, we will study now the operator whose symbol is $m\cdot C'\cdot C''$. 
It can be written as

$$T_{m\cdot C'\cdot C''}(f_1, f_2)(x)=$$

$$\int_{\R^6} m(\xi, \eta) D_{2^{k'}}^{\infty}\widehat{\phi}(\xi_1) D_{2^{k'}}^{\infty}\widehat{\psi}(\eta_1)
D_{2^{k''}}^{\infty}\widehat{\psi}(\xi_2) D_{2^{k''}}^{\infty}\widehat{\phi}(\eta_2)\widehat{f_1}(\xi_1, \xi_2)
\widehat{f_2}(\eta_1, \eta_2)e^{2\pi i x(\xi+\eta)}\,d\xi d\eta dk' dk''=$$

$$\int_{\R^6}m(\xi, \eta)\widehat{\Phi_{1, k', k''}}(\xi_1, \xi_2)\widehat{\Phi_{2, k', k''}}(\eta_1, \eta_2)
\widehat{f_1}(\xi_1, \xi_2)
\widehat{f_2}(\eta_1, \eta_2)e^{2\pi i x(\xi+\eta)}\,d\xi d\eta dk' dk''=$$

$$\int_{\R^6}m(\xi, \eta)\widehat{f_1\ast\Phi_{1, k', k''}}(\xi)\widehat{f_2\ast\Phi_{2, k', k''}}(\eta)
e^{2\pi i x(\xi+\eta)}\,d\xi d\eta dk' dk''$$
where $\Phi_{1, k', k''}:=D^1_{2^{-k'}}\phi\otimes D^1_{2^{-k''}}\psi$ and
$\Phi_{2, k', k''}:=D^1_{2^{-k'}}\psi\otimes D^1_{2^{-k''}}\phi$.

In particular, the trilinear form $\Lambda_{m\cdot C'\cdot C''}(f_1, f_2, f_3):=
\int_{\R^2}T_{m\cdot C'\cdot C''}(f_1, f_2)(x)f_3(x)\, dx$ associated to it, can be written as

\begin{equation}\label{form}
\int_{\xi+\eta+\gamma=0} m_{k', k''}(\xi, \eta, \gamma)\widehat{f_1\ast\Phi_{1, k', k''}}(\xi)\widehat{f_2\ast\Phi_{2, k', k''}}(\eta)
\widehat{f_3\ast\Phi_{3, k', k''}}(\gamma)\,d\xi d\eta d\gamma dk' dk''
\end{equation}
where $\Phi_{3, k', k''}:=D^1_{2^{-k'}}\psi'\otimes D^1_{2^{-k''}}\psi'$ and $\psi'$ is again a Schwartz function such that
$\supp(\widehat{\psi'})\subseteq [-7/4,-1/4]$ and $\widehat{\psi'}=1$ on $[-3/2,-1/2]$, while 
$m_{k', k''}(\xi, \eta, \gamma)=m(\xi, \eta)\cdot\lambda_{k', k''}(\xi, \eta, \gamma)$ where 
$\lambda_{k', k''}(\xi, \eta, \gamma)$ is a smooth function supported on
$$2\supp(\widehat{\Phi_{1, k', k''}}(\xi)\widehat{\Phi_{2, k', k''}}(\eta)\widehat{\Phi_{3, k', k''}}(\gamma))$$ which equals $1$
on $\supp(\widehat{\Phi_{1, k', k''}}(\xi)\widehat{\Phi_{2, k', k''}}(\eta)\widehat{\Phi_{3, k', k''}}(\gamma))$.

Then, we write (\ref{form}) as

$$\int_{\R^{10}}m_{k', k''}^{\vee}((n'_1, n''_1), (n'_2, n''_2), (n'_3, n''_3))\prod_{j=1}^3(f_j\ast\Phi_{j, k', k''})((x', x'')-(n'_j, n''_j))\,
dn'_j dn''_j dx' dx'' dk' dk''=$$

$$\int_{\R^{10}}2^{-4k'}2^{-4k''}m_{k', k''}^{\vee}((2^{-k'}n'_1, 2^{-k''}n''_1), 
(2^{-k'}n'_2, 2^{-k''}n''_2),(2^{-k'}n'_3, 2^{-k''}n''_3))\cdot$$

$$\prod_{j=1}^3(f_j\ast\Phi_{j, k', k''})((2^{-k'}x', 2^{-k''}x'')-(2^{-k'}n'_j, 2^{-k''}n''_j))\,
dn'_j dn''_j dx' dx'' dk' dk''=$$

$$\int_{\R^{10}}2^{-3k'}2^{-3k''}m_{k', k''}^{\vee}((2^{-k'}n'_1, 2^{-k''}n''_1), 
(2^{-k'}n'_2, 2^{-k''}n''_2),(2^{-k'}n'_3, 2^{-k''}n''_3))\cdot$$

$$2^{k'/2}2^{k''/2}\prod_{j=1}^3\langle f_j, \Phi_{j, \vec{k}, \vec{x}, \vec{n_j}}\rangle \, d\vec{n_j} d\vec{x} d\vec{k}$$
where we denoted

$$\Phi_{j, \vec{k}, \vec{x}, \vec{n_j}}:= 2^{-k'/2}2^{-k''/2}\tau_{(2^{-k'}x', 2^{-k''}x'')-(2^{-k'}n'_j, 2^{-k''}n''_j)}\Phi_{j, k', k''}.$$
Notice that our functions $\Phi_{j, \vec{k}, \vec{x}, \vec{n_j}}$ are now $L^2(\R^2)$ normalized. The above expression
can be discretized as

\begin{equation}
\sum_{(\vec{n_1}, \vec{n_2}, \vec{n_3}, \vec{k}, \vec{l})\in\Z^{10}} \Lambda_{\vec{n_1}, \vec{n_2}, \vec{n_3}, \vec{k}, \vec{l}}(f_1, f_2, f_3)
\end{equation}
where

$$\Lambda_{\vec{n_1}, \vec{n_2}, \vec{n_3}, \vec{k}, \vec{l}}(f_1, f_2, f_3):=\int_{[0,1]^{10}}2^{-3(k'+\kappa')}2^{-3(k''+\kappa'')}
m_{k'+\kappa', k''+\kappa''}^{\vee}(\cdots)$$

$$2^{(k'+\kappa')/2}2^{(k''+\kappa'')/2}\prod_{j=1}^3\langle f_j, \Phi_{j, \vec{k}+\vec{\kappa}, \vec{l}+\vec{\lambda}, \vec{n_j}+\vec{\nu_j}}\rangle \, 
d\vec{\nu_j}  d\vec{\kappa} d\vec{\lambda}.$$
Consequently, the operator $T_{m\cdot C'\cdot C''}(f_1, f_2)$ splits as

$$T_{m\cdot C'\cdot C''}(f_1, f_2)=\sum_{(\vec{n_1}, \vec{n_2}, \vec{n_3}, \vec{k}, \vec{l})\in\Z^{10}} T_{\vec{n_1}, \vec{n_2}, \vec{n_3}, \vec{k}, \vec{l}}(f_1, f_2)$$
where $T_{\vec{n_1}, \vec{n_2}, \vec{n_3}, \vec{k}, \vec{l}}$ is the operator whose trilinear form is $\Lambda_{\vec{n_1}, \vec{n_2}, \vec{n_3}, \vec{k}, \vec{l}}$. 
Clearly, by Fatou's theorem
it is enough to prove estimates for the operator

\begin{equation}\label{op1}
\sum_{(\vec{n_1}, \vec{n_2}, \vec{n_3}) \in\Z^6; |\vec{k}|, |\vec{l}|< N} T_{\vec{n_1}, \vec{n_2}, \vec{n_3}, \vec{k}, \vec{l}}(f_1, f_2)
\end{equation}
as long as they are independent of the constant $N$. Now fix a large constant $N$ and write (\ref{op1}) as

\begin{equation}\label{op2}
\sum_{(\vec{n_1}, \vec{n_2}, \vec{n_3}) \in\Z^6}\left(\sum_{|\vec{k}|, |\vec{l}|< N} T_{\vec{n_1}, \vec{n_2}, \vec{n_3}, \vec{k}, \vec{l}}(f_1, f_2)\right).
\end{equation}
We also observe that by using (\ref{2par}) and integrating by parts several times, we have
\begin{equation}\label{decay}
\left|2^{-3k'}2^{-3k''}m_{k', k''}^{\vee}((2^{-k'}n'_1, 2^{-k''}n''_1), 
(2^{-k'}n'_2, 2^{-k''}n''_2),(2^{-k'}n'_3, 2^{-k''}n''_3))\right|\lesssim
\end{equation}

$$\prod_{j=1}^3\frac{1}{(1+|\vec{n_j}|)^M},$$
for $M$ arbitrarily large.

We are going to prove explicitly that the operator
\begin{equation}\label{op3}
\sum_{|\vec{k}|, |\vec{l}|< N} T_{\vec{0}, \vec{0}, \vec{0}, \vec{k}, \vec{l}}(f_1, f_2):= \sum_{|\vec{k}|, |\vec{l}|< N} T_{\vec{k}, \vec{l}}(f_1, f_2)
\end{equation}
satisfies the required estimates. It will be clear from the proof and (\ref{decay}) that the same arguments give

\begin{equation}\label{est1}
\|\sum_{|\vec{k}|, |\vec{l}|< N} T_{\vec{n_1}, \vec{n_2}, \vec{n_3}, \vec{k}, \vec{l}}\|_{L^p\times L^q\rightarrow L^r}\lesssim
\prod_{j=1}^3\frac{1}{(1+|\vec{n_j}|)^{100}}
\|\sum_{|\vec{k}|, |\vec{l}|< N} T_{\vec{k}, \vec{l}}\|_{L^p\times L^q\rightarrow L^r}
\end{equation}
for any $(\vec{n_1}, \vec{n_2}, \vec{n_3})\in\Z^6$. Together with (\ref{op2}) 
this would prove our desired estimates. 
It is therefore enough to deal with

$$\sum_{|\vec{k}|, |\vec{l}|< N} T_{\vec{k}, \vec{l}}(f_1, f_2).$$
Fix now $p, q$ two numbers bigger than $1$ and very close to $1$. Let also $f_1, f_2$ such that $\|f_1\|_p=
\|f_2\|_q=1$. We will show that

\begin{equation}\label{est2}
\|\sum_{|\vec{k}|, |\vec{l}|< N} T_{\vec{k}, \vec{l}}(f_1, f_2)\|_{r, \infty}\lesssim 1
\end{equation}
where $1/r=1/p+1/q$.

Using Lemma 5.4 in \cite{psctc} and scaling invariance, it is enough to show that for every set $E_3\subseteq \R^2$, $|E_3|=1$,
one can find a subset $E'_3\subseteq E_3$ with $|E'_3|\sim 1$ and such that

\begin{equation}\label{est3}
|\sum_{|\vec{k}|, |\vec{l}|< N}\Lambda_{\vec{k}, \vec{l}}(f_1, f_2, f_3)|\lesssim 1
\end{equation}
where $f_3:=\chi_{E'_3}$. If this is true, then by using the symmetry of our form, the symmetry of our arguments plus multilinear interpolation as in \cite{cct1},
we would complete the proof.

In order to construct the set $E'_3$ we need to define the ''maximal-square function'' and the ''square-maximal function'' as follows.

For $(x',x'')\in\R^2$ define
$$MS(f_1)(x',x''):=\sup_{k', l'}
\frac{1}{2^{-k'_1/2}}\left(
\sum_{k'', l''}\sup_{\vec{\kappa}, \vec{\lambda}, \vec{\nu_1}}\frac{|\langle f_1, \Phi_{1, \vec{k}+\vec{\kappa}, \vec{l}+\vec{\lambda}, \vec{\nu_1}}\rangle|^2}
{2^{-k''}} 1_{I_{k'', l''}}(x'')\right)^{1/2} 1_{I_{k', l'}}(x'),$$
and

$$SM(f_2)(x',x''):=\left(
\sum_{k', k''}\frac{(\sup_{k'', l''} \sup_{\vec{\kappa}, \vec{\lambda}, \vec{\nu_2}}
\frac{|\langle f_2, \Phi_{2, \vec{k}+\vec{\kappa}, \vec{l}+\vec{\lambda}, \vec{\nu_2}}\rangle|}
{2^{-k''/2}} 1_{I_{k'', l''}}(x''))^2 }{2^{-k'}} 1_{I_{k', l'}}(x')\right)^{1/2}.$$
Then, we also define the following ''double square function''

$$SS(f_3)(x',x''):=\left(
\sum_{k', l', k'', l''}\sup_{\vec{\kappa}, \vec{\lambda}, \vec{\nu_3}}
\frac{|\langle f_3, \Phi_{3, \vec{k}+\vec{\kappa}, \vec{l}+\vec{\lambda}, \vec{\nu_3}}\rangle|^2}{2^{-k'}2^{-k''}} 1_{I_{k', l'}}(x')
1_{I_{k'', l''}}(x'')\right)^{1/2}$$
where in general $I_{k, l}$ is the dyadic interval $2^{-k}[l, l+1]$. Finally, we recall the biparameter Hardy-Littlewood
maximal function

$$MM(g)(x',x''):=\sup_{(x',x'')\in R}\frac{1}{|R|}\int_R |g(y',y'')|\, dy'dy''$$
where $R$ ranges over all rectangles in the plane, whose sides are parallel to the coordinate axes.

The reader should not worry too much about the presence of the suprema over $\kappa, \lambda,\nu_1, \nu_2, \nu_3$ in the above definitions. 
They need to be there for some technical reasons, but their appearance is 
completely harmless from the point of view of the boundedness of the corresponding operators.

It is well known that both the biparameter maximal function $MM$ and the double square function $SS$ map $L^p(\R^2)$
into $L^p(\R^2)$ whenever $1<p<\infty$, see \cite{cf}.

Similarly, it is not difficult to observe, by using Fubini's theorem and the Fefferman-Stein inequality \cite{fs},
that the operators $MS$, $SM$, are also bounded on $L^p(\R^2)$ if $1< p< \infty$ (first, one treats the $SM$ function iteratively, as we said, and then one simply
observes that the $MS$ function is pointwise smaller than $SM$).

We then set

$$\Omega_0=\{ x\in\R^2 : MS(f_1)(x)>C\} \cup \{ x\in\R^2 : SM(f_2)(x)>C\}$$
$$\cup  \{x\in\R^2 : MM(f_1)(x)>C\} \cup \{x\in\R^2 : MM(f_2)(x)>C\}.$$
Also, define

\begin{equation}\label{omega}
\Omega= \{ x\in\R^2 : MM(1_{\Omega_0})(x)>\frac{1}{100} \}
\end{equation}
and finally

$$\tilde{\Omega}= \{ x\in\R^2 : MM(1_{\Omega})(x)> \frac{1}{2} \}.$$
Clearly, we have $|\tilde{\Omega}|< 1/2$, if $C$ is a big enough constant, which we fix from now on.
Then, we define $E'_3:= E_3\setminus\tilde{\Omega}=E_3\cap\tilde{\Omega}^c$ and observe that $|E'_3|\sim 1$.

Since the form $\sum_{|\vec{k}|, |\vec{l}|< N}\Lambda_{\vec{k},\vec{l}}(f_1, f_2, f_3)$ is an average of some other forms depending
on parameters $(\vec{\kappa}, \vec{\lambda}, \vec{\nu_1}, \vec{\nu_2}, \vec{\nu_3})\in [0,1]^{10}$, it is enough to prove our inequality (\ref{est3}) for each of them, 
uniformly with respect to  $(\vec{\kappa}, \vec{\lambda}, \vec{\nu_1}, \vec{\nu_2}, \vec{\nu_3})$. We will do this in the particular case when
all these parameters are zero, but the same argument works in general. In this case, we prefer to change our notation
and write the corresponding form as

\begin{equation}\label{2paraproduct}
\Lambda_{\vec{\P}}(f_1, f_2, f_3)=\int_{\R^2} \Pi_{\vec{\P}}(f_1, f_2)(x) f_3(x)\, dx=
\sum_{\vec{P}\in\vec{\P}}\frac{1}{|I_{\vec{P}}|^{1/2}}
\langle f_1, \Phi_{\vec{P_1}}\rangle
\langle f_2, \Phi_{\vec{P_2}}\rangle
\langle f_3, \Phi_{\vec{P_3}}\rangle,
\end{equation}
where the $\vec{P}$'s are biparameter tiles corresponding to the indices $k', l', k'', l''$. More precisely, we have

$$\vec{P_1}=(P'_1, P''_1)=(2^{-k'}[l', l'+1]\times 2^{k'}[-1/4, 1/4], 2^{-k''}[l'', l''+1]\times 2^{k''}[3/4, 5/4])$$

$$\vec{P_2}=(P'_2, P''_2)=(2^{-k'}[l', l'+1]\times 2^{k'}[3/4, 5/4], 2^{-k''}[l'', l''+1]\times 2^{k''}[-1/4, 1/4])$$

$$\vec{P_3}=(P'_3, P''_3)=(2^{-k'}[l', l'+1]\times 2^{k'}[-7/4, -1/4], 2^{-k''}[l'', l''+1]\times 2^{k''}[-7/4,-1/4])$$
and $|I_{\vec{P}}|:= |I_{\vec{P_1}}|=|I_{\vec{P_2}}|=|I_{\vec{P3}}|=2^{-k'} 2^{-k''}$. 

$\P$ will be a finite set
of such biparameters tiles.
Note that $\vec{P_1}$, $\vec{P_2}$, $\vec{P_3}$ are the biparameter Heisenberg boxes of the $L^2$ normalized 
wave packets $\Phi_{\vec{P_1}}$, $\Phi_{\vec{P_2}}$, $\Phi_{\vec{P_3}}$ respectively. These new functions $\Phi_{\vec{P_j}}$ are just the old
functions $\Phi_{j,\vec{k},\vec{l}}$ previously defined, for $j=1,2,3$. We therefore need to show the following inequality

\begin{equation}\label{est4}
\sum_{\vec{P}\in\vec{\P}}\frac{1}{|I_{\vec{P}}|^{1/2}}
|\langle f_1, \Phi_{\vec{P_1}}\rangle|
|\langle f_2, \Phi_{\vec{P_2}}\rangle|
|\langle f_3, \Phi_{\vec{P_3}}\rangle|
\lesssim 1,
\end{equation}
in order to finish the proof. This will be our main goal in the next sections.

At the end of this section we would like to observe that it is very easy to obtain the desired estimates when all the indices are strictly between $1$ and $\infty$. 
To see this, let $f_1\in L^p$, $f_2\in L^q$, $f_3\in L^r$ where 
$1< p, q, r< \infty$ with $1/p+1/q+1/r=1$. Then,

$$\left|\int_{\R^2} \Pi_{\vec{\P}}(f_1, f_2)(x) f_3(x)\, dx \right|\lesssim
\sum_{\vec{P}\in\vec{\P}}\frac{1}{|I_{\vec{P}}|^{1/2}}
|\langle f_1, \Phi_{\vec{P_1}}\rangle|
|\langle f_2, \Phi_{\vec{P_2}}\rangle|
|\langle f_3, \Phi_{\vec{P_3}}\rangle|=$$

$$\int_{\R^2}\sum_{\vec{P}\in\vec{\P}}
\frac{|\langle f_1, \Phi_{\vec{P_1}}\rangle|}{|I_{\vec{P}}|^{1/2}}
\frac{|\langle f_2, \Phi_{\vec{P_2}}\rangle|}{|I_{\vec{P}}|^{1/2}}
\frac{|\langle f_3, \Phi_{\vec{P_3}}\rangle|}{|I_{\vec{P}}|^{1/2}}
\chi_{I_{\vec{P}}}(x)\, dx\lesssim$$

$$\int_{\R^2} MS(f_1)(x) SM(f_2)(x) SS(f_3)(x)\, dx\lesssim
\|MS(f_1)\|_p
\|SM(f_2)\|_q
\|SS(f_3)\|_r\lesssim$$

$$\|f_1\|_p  \|f_2\|_q \|f_3\|_r.$$

\section{Proof of the one-parameter case}

In the particular case when $\vec{\P}= \P'\times\P''$ and all the funtions $f_j$ are tensor product type functions
(i.e. $f_j= f'_j\otimes f''_j$, $j=1,2,3$), our biparameter paraproduct splits as

$$\Lambda_{\vec{\P}}(f_1, f_2, f_3)= \Lambda_{\P'}( f'_1, f'_2, f'_3) \Lambda_{\P''}(f''_1, f''_2, f''_3).$$
In this section we describe an argument which proves $L^p$ estimates for these one-parameter paraproducts
$\Lambda_{\P'}$ and $\Lambda_{\P''}$. One one hand, this method will be very useful for us in Section 6 and on the other hand it provides a new proof of the 
classical Coifman-Meyer theorem. A sketch of it in a simplified " Walsh framework " has been presented in the expository paper
\cite{psctc}.

If $I$ is an interval on the real line, we denote by $\tilde{\chi_I}(x)$ the function

$$\tilde{\chi_I}(x) = (1+ \frac{\dist(x, I)}{|I|} )^{-M}$$
where $M>0$ is a big and fixed constant. For simplicity of notation we will suppress the ''primes'' and write (for instance) $\Lambda_{\P'}(f'_1, f'_2, f'_3)$ simply as

\begin{equation}\label{paraproduct}
\Lambda_{\P}(f_1, f_2, f_3)=
\sum_{P\in\P}\frac{1}{|I_{P}|^{1/2}}
\langle f_1, \Phi_{P_1}\rangle
\langle f_2, \Phi_{P_2}\rangle
\langle f_3, \Phi_{P_3}\rangle.
\end{equation}
Notice that in this case, as $P$ runs inside the finite set $\P$, the frequency supports $\supp \widehat{\Phi_{P_j}}$,
$j=2,3$ lie inside some intervals which are essentially lacunarily disjoint, while the frequency intervals
$\supp \widehat{\Phi_{P_1}}$ are all intersecting each other.

In order to deal with this expression (\ref{paraproduct}) we need to introduce some definitions.

\begin{definition}
Let $\P$ be a finite set of tiles as before. For $j=1$ we define
$$\size_{\P}(f_j):= \sup_{P\in\P}\frac{|\langle f_j, \Phi_{P_j}\rangle|}{|I_P|^{1/2}}$$
and for $j=2,3$ we set

$$\size_{\P}(f_j):= \sup_{P\in\P}\frac{1}{|I_P|}
\left\|\left(\sum_{I_{P'}\subseteq I_P}\frac{|\langle f_j, \Phi_{P'_j}\rangle |^2}{|I_{P'}|} 1_{I_{P'}}\right)^{1/2}\right\|_{1,\infty}.$$
Also, for $j=1,2,3$, we define

$$\energy_{\P}(f_j):=\sup_{\D\subseteq\P}
\left\|\sum_{P\in \D}\frac{\langle |f_j|, \tilde{\chi_{I_P}}\rangle}{|I_P|} 1_{I_P}\right\|_{1,\infty}$$
where $\D$ ranges over all subsets of $\P$ such that the intervals $\{I_P : P\in\D\}$ are disjoint.
\end{definition}

The following John-Nirenberg type inequality holds in this context (see \cite{cct1}).

\begin{lemma}\label{jn}
Let $\P$ be a finite collection of tiles as before and $j=2,3$. Then

$$\size_{\P}(f_j) \sim \sup_{P\in\P}
\left(\frac{1}{|I_P|}\sum_{I_{P'}\subseteq I_P}|\langle f_j, \Phi_{P'_j}\rangle|^2\right)^{1/2}.$$
\end{lemma}
We will also need the following lemma (see also \cite{cct1}).

\begin{lemma}\label{cz}
Let $\P$ be a finite collection of tiles and $j=2,3$. Then, we have

$$\left\|\left(\sum_{I_{P'}\subseteq I_P}\frac{|\langle f, \Phi_{P'_j}\rangle |^2}{|I_{P'}|} 1_{I_{P'}}\right)^{1/2}\right\|_{1,\infty}\lesssim \|f\tilde{\chi_{I_P}}\|_1.$$
\end{lemma}
The following proposition will be very helpful.

\begin{proposition}
Let $j=1, 2, 3$, $\P'$ a subset of $\P$, $n\in \Z$ and suppose that

$$\size_{\P'}(f_j)\leq 2^{-n} \energy_{\P}(f_j).$$
Then, we may decompose $\P'= \P'' \cup \P'''$ such that

\begin{equation}\label{size}
\size_{\P''}(f_j)\leq 2^{-n-1} \energy_{\P}(f_j)
\end{equation}
and that $\P'''$ can be written as a disjoint union of subsets $T\in \T$ such that for every $T\in \T$, there exists
an interval $I_T$ (corresponding to a certain tile) having the property that every $P\in T$ has $I_P\subseteq I_T$
and also such that

\begin{equation}\label{L1}
\sum_{T\in \T} |I_T| \lesssim 2^n.
\end{equation}
\end{proposition}

\begin{proof}
The idea is to remove large subsets of $\P'$ one by one, placing them into $\P'''$ until (\ref{size}) is satisfied.

\underline{Case 1}: $j=1$. Pick a tile $P\in \P'$ such that $|I_P|$ is as big as possible and such that

$$\frac{|\langle f_j, \Phi_{P_j}\rangle|}{|I_P|^{1/2}} > 2^{-n-1} \energy_{\P}(f_j).$$
Then, collect all the tiles $P'\in\P'$ such that $I_{P'}\subseteq I_P$ into a set called $T$ and place $T$ into
$\P'''$. Define $I_T:= I_P$. Then look at the remaining tiles in $\P'\setminus T$ and repeat the procedure. Since
there are finitely many tiles, the procedure ends after finitely many steps producing the subsets $T\in \T$.
Clealy, (\ref{size}) is now satisfied and it remains to show (\ref{L1}). To see this, one can write

$$\sum_{T\in \T} |I_T| = \|\sum_{T\in \T} 1_{I_T} \|_1 =
\|\sum_{T\in \T} 1_{I_T} \|_{1,\infty}$$
since by construction, our intervals $I_T$ are disjoint. Then, the right hand side of the above equality is smaller than

$$2^n \energy_{\P}(f_j)^{-1} \left\|\sum_{T\in \T}\frac{\langle |f_j|, \tilde{\chi_{I_T}}\rangle}{|I_T|} 1_{I_T}\right\|_{1,\infty}\lesssim 2^n.$$

\underline{Case 2}: $j=2,3$. The algorithm is very similar. Pick again a tile $P\in \P'$ such that
$|I_P|$ is as big as possible and such that

$$\frac{1}{|I_P|}
\left\|\left(\sum_{I_{P'}\subseteq I_P}\frac{|\langle f_j, \Phi_{P'_j}\rangle |^2}{|I_{P'}|} 1_{I_{P'}}\right)^{1/2}\right\|_{1,\infty} > 2^{-n-1} \energy_{\P}(f_j).$$
Then, as before, collect all the tiles $P'\in \P'$ such that $I_{P'}\subseteq I_P$ in a set named $T$and place this $T$
into $\P'''$. Define, as in Case 1, $I_T:= I_P$. Then look at the remaining tiles $\P'\setminus T$ and repeat the procedure which of course ends after finitely many steps. Inequality (\ref{size}) is now clear, it remains to understand
(\ref{L1}) only.

Since the intervals $I_T$ are disjoint by construction, we can write

$$\sum_{T\in \T} |I_T| = \|\sum_{T\in \T} 1_{I_T} \|_1 =
\|\sum_{T\in \T} 1_{I_T} \|_{1,\infty}\lesssim$$

$$2^n \energy_{\P}(f_j)^{-1}
\left\|\sum_{T\in \T}
\frac{1}{|I_T|}
\|(\sum_{I_{P'}\subseteq I_T}\frac{|\langle f_j, \Phi_{P'_j}\rangle |^2}{|I_{P'}|} 1_{I_{P'}})^{1/2}\|_{1,\infty} 1_{I_T}\right\|_{1,\infty}\lesssim$$

$$2^n \energy_{\P}(f_j)^{-1}\left\|
\sum_{T\in \T}\frac{\langle |f_j|, \tilde{\chi_{I_T}}\rangle}{|I_T|} 1_{I_T} \right\|_{1, \infty}\lesssim 2^n,$$
by using Lemma \ref{cz}, and this ends the proof.
\end{proof}

By iterating the above lemma, we immediately obtain the following consequence.

\begin{corollary}\label{use}
Let $j=1, 2, 3$. There exists a partition

$$\P= \bigcup_{n\in \Z} \P_n$$
such that for every $n\in \Z$ we have

$$\size_{\P_n}(f_j)\leq \min (2^{-n} \energy_{\P}(f_j), \size_{\P}(f_j) ).$$
Also, we may write each $\P_n$ as a disjoint union of subsets $T\in\T_n$ as before, such that

$$\sum_{T\in\T_n} |I_T|\lesssim 2^n.$$
\end{corollary}
We now prove the following proposition.

\begin{proposition}\label{USE}
Let $\P$ be a set as before. Then,

\begin{equation}\label{in}
\sum_{P\in\P}\frac{1}{|I_{P}|^{1/2}}
|\langle f_1, \Phi_{P_1}\rangle|
|\langle f_2, \Phi_{P_2}\rangle|
|\langle f_3, \Phi_{P_3}\rangle|\lesssim \prod_{j=1}^3 \size_{\P}(f_j)^{1-\theta_j} \energy_{\P}(f_j)^{\theta_j}
\end{equation}
for any $0\leq \theta_1, \theta_2, \theta_3 <1$ such that $\theta_1+ \theta_2+ \theta_3= 1$ with the implicit
constant depending on $\theta_j$, $j=1,2,3.$
\end{proposition}

\begin{proof}
During this proof, we will write for simplicity $S_j:= \size_{\P}(f_j)$ and $E_j:= \energy_{\P}(f_j)$, for
$j=1,2,3$. If we apply Corollary \ref{use} to the functions $\frac{f_j}{E_j}$, $j=1,2,3$ we obtain
a decomposition
 
$$\P=\bigcup_n \P^j_n$$
such that each $\P^j_n$ can be written as a union of subsets in $\T^j_n$ with the properties described in 
Corollary \ref{use}.
In particular, one can write the left hand side of our described inequality (\ref{in}) as

\begin{equation}\label{in1}
E_1 E_2 E_3 \sum_{n_1, n_2, n_3} \sum_{T\in \T^{n_1, n_2, n_3}} 
\sum_{P\in T}\frac{1}{|I_{P}|^{1/2}}
|\langle \frac{f_1}{E_1}, \Phi_{P_1}\rangle|
|\langle \frac{f_2}{E_2}, \Phi_{P_2}\rangle|
|\langle \frac{f_3}{E_3}, \Phi_{P_3}\rangle|
\end{equation}
where $\T^{n_1, n_2, n_3}:= \T^1_{n_1} \cap \T^2_{n_2} \cap \T^3_{n_3}$.  By using H\"{o}lder inequality on every $T\in\T^{n_1, n_2, n_3}$ together with Lemma \ref{jn}, one can estimate the sum in (\ref{in1}) by

\begin{equation}\label{in2}
E_1 E_2 E_3\sum_{n_1, n_2, n_3} 2^{-n_1} 2^{-n_2} 2^{-n_3} \sum_{T\in\T^{n_1, n_2, n_3}}|I_T|
\end{equation}
where (according to the same Corollary \ref{use}) the summation goes over those $n_1, n_2, n_3 \in \Z$ satisfying

\begin{equation}\label{must}
2^{-n_j}\lesssim \frac{S_j}{E_j}.
\end{equation}
On the other hand, Corollary \ref{use} allows us to estimate the inner sum in (\ref{in2}) in three different ways,
namely

$$\sum_{T\in\T^{n_1, n_2, n_3}} |I_T|\lesssim 2^{n_1}, 2^{n_2}, 2^{n_3}$$
and so, in particular, we can also write

\begin{equation}\label{in3}
\sum_{T\in\T^{n_1, n_2, n_3}} |I_T|\lesssim 2^{n_1\theta_1} 2^{n_2\theta_2} 2^{n_3\theta_3}
\end{equation}
whenever $0\leq \theta_1, \theta_2, \theta_3 <1$ with $\theta_1+ \theta_2+ \theta_3= 1$. Using (\ref{in3}) and (\ref{must}),
one can estimate (\ref{in2}) further by

$$E_1 E_2 E_3 \sum_{n_1, n_2, n_3}2^{-n_1(1-\theta_1)} 2^{-n_2(1-\theta_2)} 2^{-n_3(1-\theta_3)}\lesssim$$

$$E_1 E_2 E_3 
(\frac{S_1}{E_1})^{1-\theta_1}
(\frac{S_2}{E_2})^{1-\theta_2} 
(\frac{S_2}{E_2})^{1-\theta_3}= \prod_{j=1}^3 S_j^{1-\theta_j} \prod_{j=1}^3 E_j^{\theta_j},$$
which ends the proof.

\end{proof}

Using this Proposition \ref{USE}, one can prove the $L^p$ boundedness of one-parameter paraproducts, as follows.
We just need to show that they map $L^1\times L^1\rightarrow L^{1/2, \infty}$, because then, by interpolation and symmetry one can deduce that they map $L^p\times L^q\rightarrow L^r$ as long as $1<p,q \leq \infty$, $0<r<\infty$ and
$1/p+1/q=1/r$.

Let $f_1, f_2 \in L^1$ be such that $\|f_1\|_1= \|f_2\|_1=1$. As before, it is enough to show that given
$E_3\subseteq \R$ $|E_3|=1$, one can find a subset $E'_3\subseteq E_3$ with $|E'_3|\sim 1$ and

\begin{equation}\label{in4}
\sum_{P\in\P}\frac{1}{|I_{P}|^{1/2}}
|\langle f_1, \Phi_{P_1}\rangle|
|\langle f_2, \Phi_{P_2}\rangle|
|\langle f_3, \Phi_{P_3}\rangle|\lesssim 1
\end{equation}
where $f_3:=\chi_{E'_3}$. For, we define the set $U$ by

$$U:= \{ x\in\R : M(f_1)(x)> C\} \cup \{ x\in \R : M(f_2)(x)> C \}$$
where $M(f)$ is the Hardy-Littlewood maximal operator of $f$. Clearly, we have $|U|< 1/2$ if $C>0$ is big enough. Then we
define our set $E'_3:=E_3\cap U^c$ and remark that $|E'_3|\sim 1$.

Then, we write 

$$\P= \bigcup_{d\geq 0} \P_d,$$
where

$$\P_d:= \{ P\in \P : \frac{\dist(I_P, U^c)}{|I_P|}\sim 2^d \}.$$
After that, by using Lemma \ref{cz}, we observe that $\size_{\P_d}(f_j)\lesssim 2^d$ for $j=1,2$, while
$\size_{\P_d}(f_3)\lesssim 2^{-Nd}$ for an arbirarily big number $N>0$. We also observe that

$$\energy_{\P_d}(f_j)\lesssim \|M(f_j)\|_{1,\infty} \lesssim \|f_j\|_1 = 1.$$
By applying Proposition \ref{USE} in the particular case $\theta_1= \theta_2= \theta_3=1/3$, we get
that the left hand side of (\ref{in4}) can be majorized by

$$\sum_{d\geq 0}\sum_{P\in \P_d}
\frac{1}{|I_{P}|^{1/2}}
|\langle f_1, \Phi_{P_1}\rangle|
|\langle f_2, \Phi_{P_2}\rangle|
|\langle f_3, \Phi_{P_3}\rangle|\lesssim \sum_{d\geq 0} 2^{2d/3} 2^{2d/3} 2^{-2Nd/3}\lesssim 1$$
as wanted and this finishes the proof of the one-parameter case.

The reader should compare this Proposition \ref{USE} with the corresponding Proposition 6.5 in \cite{cct3}. Our present ``lacunary setting'' allows
for an $L^1$-type definition of the ``energies'' (instead of $L^2$-type as in \cite{cct3}) and this is why we can obtain the full range of estimates this time.

\section{Proof of Theorem \ref{main}}
  
We reduced our proof to showing (\ref{est4}). Clearly, this inequality is the bi-parameter analogue of the inequality (\ref{in4}) above. Unfortunately, the technique just
described in Section 3, so useful when estimating (\ref{in4}), cannot handle our sum in (\ref{est4}) this time. In fact, we do not know if there exists a satisfactory 
bi-parameter analogue of Proposition \ref{USE} 
and this is where some of the main new difficulties are coming from. Hence, we have to proceed differently.

We split the left hand side of that inequality
into two parts, as follows

\begin{equation}\label{III}
\sum_{\vec{P}}= \sum_{I_{\vec{P}}\cap\Omega^c\neq\emptyset} + \sum_{I_{\vec{P}}\cap\Omega^c = \emptyset}
:= I + II
\end{equation}
where $\Omega$ is the set defined in (\ref{omega}).
\section{Estimates for term I}

We first estimate term $I$. 
The argument goes as follows.

Since $I_{\vec{P}}\cap\Omega^c\neq\emptyset$, it follows that 
$\frac{|I_{\vec{P}}\cap\Omega_0|}{|I_{\vec{P}}|}< \frac{1}{100}$ or equivalently,
$|I_{\vec{P}}\cap\Omega_0^c|> \frac{99}{100}|I_{\vec{P}}|$. 

We are now going to describe three decomposition procedures, one for each function $f_1, f_2, f_3$. Later on, we will
combine them, in order to handle our sum.

First, define 

$$\Omega_1= \{ x\in\R^2 : MS(f_1)(x)> \frac{C}{2^1} \}$$
and set

$$\T_1= \{ \vec{P}\in \vec{\P} : |I_{\vec{P}}\cap\Omega_1|>\frac{1}{100} |I_{\vec{P}}| \},$$
then define

$$\Omega_2= \{ x\in\R^2 : MS(f_1)(x)> \frac{C}{2^2} \}$$
and set

$$\T_2= \{ \vec{P}\in \vec{\P}\setminus\T_1 : |I_{\vec{P}}\cap\Omega_2|>\frac{1}{100} |I_{\vec{P}}| \},$$
and so on. The constant $C>0$ is the one in the definition of the set $E'_3$ in Section 2.
Since there are finitely many tiles, this algorithm ends after a while, producing the sets $\{\Omega_n\}$
and $\{\T_n\}$ such that $\vec{\P}=\cup_n\T_n$.

Independently, define

$$\Omega'_1= \{ x\in\R^2 : SM(f_2)(x)> \frac{C}{2^1} \}$$
and set

$$\T'_1= \{ \vec{P}\in \vec{\P} : |I_{\vec{P}}\cap\Omega'_1|>\frac{1}{100} |I_{\vec{P}}| \},$$
then define

$$\Omega'_2= \{ x\in\R^2 : SM(f_2)(x)> \frac{C}{2^2} \}$$
and set

$$\T'_2= \{ \vec{P}\in \vec{\P}\setminus\T'_1 : |I_{\vec{P}}\cap\Omega'_2|>\frac{1}{100} |I_{\vec{P}}| \},$$
and so on, producing the sets $\{\Omega'_n\}$ and $\{\T'_n\}$ such that $\vec{\P}=\cup_n\T'_n$.
We would like to have such a decomposition available for the function $f_3$ also. To do this, we first need to
construct the analogue of the set $\Omega_0$, for it. Pick  $N>0$ a big enough integer such that for every
$\vec{P}\in\vec{\P}$ we have $|I_{\vec{P}}\cap\Omega^{''c}_{-N}|> \frac{99}{100} |I_{\vec{P}}|$ where we defined

$$\Omega''_{-N}= \{ x\in\R^2 : SS(f_3)(x)> C 2^N \}.$$
Then, similarly to the previous algorithms, we define

$$\Omega''_{-N+1}= \{ x\in\R^2 : SS(f_2)(x)> \frac{C 2^N}{2^1} \}$$
and set

$$\T''_{-N+1}= \{ \vec{P}\in \vec{\P} : |I_{\vec{P}}\cap\Omega''_{-N+1}|>\frac{1}{100} |I_{\vec{P}}| \},$$
then define

$$\Omega''_{-N+2}= \{ x\in\R^2 : SS(f_3)(x)> \frac{C 2^N}{2^2} \}$$
and set

$$\T''_{-N+2}= \{ \vec{P}\in \vec{\P}\setminus\T''_{-N+1} : |I_{\vec{P}}\cap\Omega''_{-N+2}|>\frac{1}{100} |I_{\vec{P}}| \},$$
and so on, constructing the sets $\{\Omega''_n\}$ and $\{\T''_n\}$ such that $\vec{\P}=\cup_n\T''_n$.

Then we write the term $I$ as

\begin{equation}\label{in5}
\sum_{n_1, n_2>0, n_3>-N} \sum_{\vec{P}\in \T_{n_1, n_2, n_3}}
\frac{1}{|I_{\vec{P}}|^{3/2}}
|\langle f_1, \Phi_{\vec{P_1}}\rangle|
|\langle f_2, \Phi_{\vec{P_2}}\rangle|
|\langle f_3, \Phi_{\vec{P_3}}\rangle| |I_{\vec{P}}|,
\end{equation}
where $\T_{n_1, n_2, n_3}:= \T_{n_1}\cap \T'_{n_2}\cap \T''_{n_3}$. Now, if $\vec{P}$ belongs to
$\T_{n_1, n_2, n_3}$ this means in particular that $\vec{P}$ has not been selected at the previous $n_1 -1$, $n_2 -1$ and
$n_3 -1$ steps respectively, which means that $|I_{\vec{P}}\cap\Omega_{n_1-1}|<\frac{1}{100} |I_{\vec{P}}|$,
$|I_{\vec{P}}\cap\Omega'_{n_2-1}|<\frac{1}{100} |I_{\vec{P}}|$ and $|I_{\vec{P}}\cap\Omega''_{n_3-1}|<\frac{1}{100} |I_{\vec{P}}|$ or equivalently, 
$|I_{\vec{P}}\cap\Omega^c_{n_1-1}|>\frac{99}{100} |I_{\vec{P}}|$,
$|I_{\vec{P}}\cap\Omega^{'c}_{n_2-1}|>\frac{99}{100} |I_{\vec{P}}|$ and
$|I_{\vec{P}}\cap\Omega^{''c}_{n_3-1}|>\frac{99}{100} |I_{\vec{P}}|$. But this implies that

\begin{equation}\label{in6}
|I_{\vec{P}}\cap\Omega^c_{n_1-1}\cap\Omega^{'c}_{n_2-1}\cap\Omega^{''c}_{n_3-1}|> \frac{97}{100}|I_{\vec{P}}|.
\end{equation}
In particular, using (\ref{in6}), the term in (\ref{in5}) is smaller than

$$\sum_{n_1, n_2>0, n_3>-N} \sum_{\vec{P}\in \T_{n_1, n_2, n_3}}
\frac{1}{|I_{\vec{P}}|^{3/2}}
|\langle f_1, \Phi_{\vec{P_1}}\rangle|
|\langle f_2, \Phi_{\vec{P_2}}\rangle|
|\langle f_3, \Phi_{\vec{P_3}}\rangle| |I_{\vec{P}}\cap\Omega^c_{n_1-1}\cap\Omega^{'c}_{n_2-1}\cap\Omega^{''c}_{n_3-1}|=$$

$$\sum_{n_1, n_2>0, n_3>-N} \int_{\Omega^c_{n_1-1}\cap\Omega^{'c}_{n_2-1}\cap\Omega^{''c}_{n_3-1}}
\sum_{\vec{P}\in \T_{n_1, n_2, n_3}}
\frac{1}{|I_{\vec{P}}|^{3/2}}
|\langle f_1, \Phi_{\vec{P_1}}\rangle|
|\langle f_2, \Phi_{\vec{P_2}}\rangle|
|\langle f_3, \Phi_{\vec{P_3}}\rangle| \chi_{I_{\vec{P}}}(x)\, dx $$

$$\lesssim \sum_{n_1, n_2>0, n_3>-N} \int_{\Omega^c_{n_1-1}\cap\Omega^{'c}_{n_2-1}\cap\Omega^{''c}_{n_3-1}\cap
\Omega_{\T_{n_1, n_2, n_3}}}
MS(f_1)(x) SM(f_2)(x) SS(f_3)(x)\, dx$$

\begin{equation}\label{in7}
\lesssim \sum_{n_1, n_2>0, n_3>-N} 2^{-n_1} 2^{-n_2} 2^{-n_3} |\Omega_{\T_{n_1, n_2, n_3}}|,
\end{equation}
where

$$\Omega_{\T_{n_1, n_2, n_3}}:= \bigcup_{\vec{P}\in\T_{n_1, n_2, n_3}} I_{\vec{P}}.$$
On the other hand we can write

$$|\Omega_{\T_{n_1, n_2, n_3}}|\leq |\Omega_{\T_{n_1}}|\leq
|\{ x\in\R^2 : MM(\chi_{\Omega_{n_1}})(x)> \frac{1}{100} \}|$$

$$\lesssim |\Omega_{n_1}|= |\{ x\in\R^2 : MS(f_1)(x)>\frac{C}{2^{n_1}} \}|\lesssim 2^{n_1 p}.$$
Similarly, we have

$$|\Omega_{\T_{n_1, n_2, n_3}}|\lesssim 2^{n_2 q}$$
and also

$$|\Omega_{\T_{n_1, n_2, n_3}}|\lesssim 2^{n_2 \alpha},$$
for every $\alpha >1$. Here we used the fact that all the operators $SM$, $MS$, $SS$, $MM$ are bounded
on $L^s$ as long as $1<s< \infty$ and also that $|E'_3|\sim 1$.
In particular, it follows that

\begin{equation}\label{*}
|\Omega_{\T_{n_1, n_2, n_3}}|\lesssim 2^{n_1 p \theta_1}
 2^{n_2 q \theta_2} 2^{n_3 \alpha \theta_3}
\end{equation}
for any $0\leq \theta_1, \theta_2, \theta_3 < 1$, such that $\theta_1+ \theta_2 +\theta_3= 1$.

Now we split the sum in (\ref{in7}) into

\begin{equation}\label{last}
\sum_{n_1, n_2>0, n_3>0} 2^{-n_1} 2^{-n_2} 2^{-n_3} |\Omega_{\T_{n_1, n_2, n_3}}|+ 
\sum_{n_1, n_2>0, 0>n_3>-N} 2^{-n_1} 2^{-n_2} 2^{-n_3} |\Omega_{\T_{n_1, n_2, n_3}}|.
\end{equation}
To estimate the first term in (\ref{last}) we use the inequality (\ref{*}) in the particular case
$\theta_1=\theta_2=1/2$, $\theta_3=0$, while to estimate the second term we use (\ref{*}) for $\theta_j$, $j=1,2,3$
such that $1-p\theta_1>0$, $1-q\theta_2>0$ and $\alpha\theta_3 -1>0$. With these choices, the sum in (\ref{last})
is $O(1)$. This ends the discussion of $I$.

\section{Estimates for term II}

It remains to estimate term $II$ in (\ref{III}). The sum now runs over those tiles having the property that
$I_{\vec{P}}\subseteq \Omega$. For every such $\vec{P}$ there exists a maximal dyadic rectangle $R$ such that
$I_{\vec{P}}\subseteq R\subseteq\Omega$. We collect all such distinct maximal rectangles into a set called
$\it{R}_{\max}$. For $d\geq 1$ an integer, we denote by $\it{R}^d_{\max}$ the set of all $R\in\it{R}_{\max}$
such that $2^d R\subseteq\tilde{\Omega}$ and $d$ is maximal with this property. 

By using Journ\'e's Lemma \cite{journe} in the form presented in \cite{fl}, we have that for every $\epsilon >0$

\begin{equation}\label{jl}
\sum_{R\in\it{R}^d_{\max}} |R|\lesssim 2^{\epsilon d} |\Omega|.
\end{equation}
Our initial sum in $II$ is now smaller than

\begin{equation}\label{sum}
\sum_{d\geq 1} \sum_{R\in\it{R}^d_{\max}}
\sum_{I_{\vec{P}}\subseteq R\cap\Omega}
\frac{1}{|I_{\vec{P}}|^{1/2}}
|\langle f_1, \Phi_{\vec{P_1}}\rangle|
|\langle f_2, \Phi_{\vec{P_2}}\rangle|
|\langle f_3, \Phi_{\vec{P_3}}\rangle|.
\end{equation}
We claim that for every $R\in\it{R}^d_{\max}$ we have

\begin{equation}\label{claim}
\sum_{I_{\vec{P}}\subseteq R\cap\Omega}
\frac{1}{|I_{\vec{P}}|^{1/2}}
|\langle f_1, \Phi_{\vec{P_1}}\rangle|
|\langle f_2, \Phi_{\vec{P_2}}\rangle|
|\langle f_3, \Phi_{\vec{P_3}}\rangle|\lesssim 2^{-N d} |R|,
\end{equation}
for any number $N>0$.
If (\ref{claim}) is true, then by combining it with (\ref{jl}), we can estimate (\ref{sum}) by

$$\sum_{d\geq 1}\sum_{R\in\it{R}^d_{\max}} 2^{-N d} |R|
= \sum_{d\geq 1}2^{-N d}\sum_{R\in\it{R}^d_{\max}} |R|$$

$$\lesssim \sum_{d\geq 1} 2^{-N d} 2^{\epsilon d}\lesssim 1,$$
which would complete the proof.

It remains to prove (\ref{claim}). Fix $R:= I\times J$ in $\it{R}^d_{\max}$. Since 
$2^d R:=\tilde{I}\times\tilde{J}\subseteq \tilde{\Omega}$,
it follows that $2^d R\cap E'_3 = \emptyset$ and so $\chi_{E'_3}=\chi_{E'_3}\chi_{(\tilde{I}\times\tilde{J})^c}$.
Now we write

$$\chi_{(\tilde{I}\times\tilde{J})^c}=\chi_{\tilde{I}^c}+\chi_{\tilde{J}^c}-\chi_{\tilde{I}^c}\cdot\chi_{\tilde{J}^c}.$$
As a consequence, the left hand side in (\ref{claim}) splits into three sums. Since all are similar, we will treat only the first one.

Recall that every $I_{\vec{P}}$ is of the form $I_{\vec{P}}= I_{P'}\times I_{P''}$ and let us denote by $\L$ the set

$$\L := \{ I_{P'} : I_{\vec{P}}\subseteq R \}.$$
Then split 

$$\L= \bigcup_{d_1\geq 0} \L_{d_1}$$
where

$$\L_{d_1}:= \{ K'\in \L : \frac{|I|}{|K'|}\sim 2^{d_1} \}$$
and observe that

\begin{equation}\label{interval}
\sum_{K'\in \L_{d_1}} |K'| \lesssim |I|.
\end{equation}
Then, we can majorize the left hand side of (\ref{claim}) by 

$$\sum_{d_1\geq 0} \sum_{K'\in \L_{d_1}} \sum_{I_{\vec{P}}\subseteq R; I_{P'}= K'}
\frac{1}{|I_{\vec{P}}|^{1/2}}
|\langle f_1, \Phi_{\vec{P_1}}\rangle|
|\langle f_2, \Phi_{\vec{P_2}}\rangle|
|\langle f_3, \Phi_{\vec{P_3}}\rangle|=$$

$$\sum_{d_1\geq 0} \sum_{K'\in \L_{d_1}} \sum_{I_{\vec{P}}\subseteq R; I_{P'}= K'}
|I_{P'}|\frac{1}{|I_{P''}|^{1/2}}
|\langle \frac{\langle f_1, \Phi_{P'_1}\rangle}{|I_{P'}|^{1/2}}, \Phi_{P''_1}\rangle|
|\langle \frac{\langle f_2, \Phi_{P'_2}\rangle}{|I_{P'}|^{1/2}}, \Phi_{P''_2}\rangle|
|\langle \frac{\langle f_3, \Phi_{P'_3}\rangle}{|I_{P'}|^{1/2}}, \Phi_{P''_3}\rangle|,$$
where we redefined $f_3:=\chi_{E'_3}\cdot\chi_{\tilde{I}^c}$.

Let us observe that if $\vec{P}$ is such that $I_{P'}=K'$ then the one-parameter tiles
$P'_j$, $j=1,2,3$ are fixed and we will denote for simplicity $\Phi_{P'_j}:=\Phi^j_{K'}$. We also denote
by

$$\P_{K'}:= \{ P'' : I_{\vec{P}}\subseteq R,\,\,\,I_{P'}=K' \}.$$
With these notations, we rewrite our sum as

\begin{equation}\label{split}
\sum_{d_1\geq 0} \sum_{K'\in \L_{d_1}}|K'| \sum_{P''\in\P_{K'}}
\frac{1}{|I_{P''}|^{1/2}}\prod_{j=1}^3
|\langle \frac{\langle f_j, \Phi^j_{K'}\rangle}{|K'|^{1/2}}, \Phi_{P''_j}\rangle|.
\end{equation} 
Next we split $\P_{K'}$ as

$$\P_{K'}= \bigcup_{d_2\geq 0}\P_{K'}^{d_2}$$
where

$$\P_{K'}^{d_2}:= \{ P''\in \P_{K'} : \frac{|J|}{|I_{P''}|}\sim 2^{d_2} \}.$$
As a consequence, (\ref{split}) splits into

\begin{equation}\label{split1}
\sum_{d_1\geq 0} \sum_{K'\in \L_{d_1}}|K'| \sum_{d_2\geq 0}\sum_{P''\in\P_{K'}^{d_2}}
\frac{1}{|I_{P''}|^{1/2}}\prod_{j=1}^3
|\langle \frac{\langle f_j, \Phi^j_{K'}\rangle}{|K'|^{1/2}}, \Phi_{P''_j}\rangle|=
\end{equation}

$$\sum_{d_1\geq 0} \sum_{K'\in \L_{d_1}}|K'| \sum_{P''\in\bigcup_{d_2\leq d_1}\P_{K'}^{d_2}}
\frac{1}{|I_{P''}|^{1/2}}\prod_{j=1}^3
|\langle \frac{\langle f_j, \Phi^j_{K'}\rangle}{|K'|^{1/2}}, \Phi_{P''_j}\rangle| +$$

$$\sum_{d_1\geq 0} \sum_{K'\in \L_{d_1}}|K'| \sum_{P''\in\bigcup_{d_2\geq d_1}\P_{K'}^{d_2}}
\frac{1}{|I_{P''}|^{1/2}}\prod_{j=1}^3
|\langle \frac{\langle f_j, \Phi^j_{K'}\rangle}{|K'|^{1/2}}, \Phi_{P''_j}\rangle|.$$
To estimate the first term on the right hand side of (\ref{split1}) we observe that

$$\size_{\bigcup_{d_2\leq d_1}\P_{K'}^{d_2}}(\frac{\langle f_1, \Phi^1_{K'} \rangle}{|K'|^{1/2}})\lesssim 2^{d_1 +d},$$
$$\size_{\bigcup_{d_2\leq d_1}\P_{K'}^{d_2}}(\frac{\langle f_2, \Phi^2_{K'} \rangle}{|K'|^{1/2}})\lesssim 2^{d_1 +d},$$
$$\size_{\bigcup_{d_2\leq d_1}\P_{K'}^{d_2}}(\frac{\langle f_3, \Phi^3_{K'} \rangle}{|K'|^{1/2}})\lesssim 2^{-N(d_1 +d)},$$
where $N$ is as big as we want. Similarly, we have

$$\energy_{\bigcup_{d_2\leq d_1}\P_{K'}^{d_2}}(\frac{\langle f_1, \Phi^1_{K'} \rangle}{|K'|^{1/2}})\lesssim 2^{d_1 +d} |J|,$$
$$\energy_{\bigcup_{d_2\leq d_1}\P_{K'}^{d_2}}(\frac{\langle f_2, \Phi^2_{K'} \rangle}{|K'|^{1/2}})\lesssim 2^{d_1 +d} |J|,$$
$$\energy_{\bigcup_{d_2\leq d_1}\P_{K'}^{d_2}}(\frac{\langle f_3, \Phi^3_{K'} \rangle}{|K'|^{1/2}})\lesssim 2^{-N(d_1 +d)} |J|.$$
Using these inequalities and applying Proposition \ref{USE},
we can majorize that first term by

\begin{equation}
\sum_{d_1\geq 0} \sum_{K'\in \L_{d_1}} |K'| 2^{d_1 +d} 2^{d_1+ d} 2^{-N(d_1+ d)} |J|=
\end{equation}

$$2^{-(N-2)d} |J|\sum_{d_1\geq 0} 2^{-(N-2)d_1}\sum_{K'\in \L_{d_1}} |K'|\lesssim$$

$$2^{-(N-2)d} |J|\sum_{d_1\geq 0} 2^{-(N-2)d_1} |I|\lesssim 2^{-(N-2)d}|I||J|=
2^{-(N-2)d} |R|,$$
also by using (\ref{interval}). Then, to handle the second term on the right hand side of (\ref{split1}), we decompose

\begin{equation}\label{d3}
\bigcup_{d_2\geq d_1}\P_{K'}^{d_2} = \bigcup_{d_3} \P_{K', d_3}
\end{equation}
where $\P_{K', d_3}$ is the collection of all tiles $P''\in \bigcup_{d_2\geq d_1}\P_{K'}^{d_2}$ so that
$2^{d_3}(K'\times I_{P''})\subseteq \tilde{\Omega}$ and $d_3$ is maximal with this property. 

It is not difficult to observe that in fact we have the constraint $d_1+d\leq d_3$.
Taking this into account, the second term can be written as

\begin{equation}\label{cinu}
\sum_{d_1\geq 0} \sum_{K'\in \L_{d_1}}|K'| \sum_{d_3\geq d_1+d}\sum_{P''\in\P_{K',d_3}}
\frac{1}{|I_{P''}|^{1/2}}\prod_{j=1}^3
|\langle \frac{\langle f_j, \Phi^j_{K'}\rangle}{|K'|^{1/2}}, \Phi_{P''_j}\rangle|.
\end{equation}

Now we estimate as before the sizes and energies as follows

$$\size_{\P_{K',d_3}}(\frac{\langle f_1, \Phi^1_{K'} \rangle}{|K'|^{1/2}})\lesssim 2^{d_3},$$
$$\size_{\P_{K',d_3}}(\frac{\langle f_2, \Phi^2_{K'} \rangle}{|K'|^{1/2}})\lesssim 2^{d_3},$$
$$\size_{\P_{K',d_3}}(\frac{\langle f_3, \Phi^3_{K'} \rangle}{|K'|^{1/2}})\lesssim 2^{-Nd_3},$$
where, as usual, $N$ is as big as we want. Similarly, we have

$$\energy_{\P_{K',d_3}}(\frac{\langle f_1, \Phi^1_{K'} \rangle}{|K'|^{1/2}})\lesssim 2^{d_3} |J|,$$
$$\energy_{\P_{K',d_3}}(\frac{\langle f_2, \Phi^2_{K'} \rangle}{|K'|^{1/2}})\lesssim 2^{d_3} |J|,$$
$$\energy_{\P_{K',d_3}}(\frac{\langle f_3, \Phi^3_{K'} \rangle}{|K'|^{1/2}})\lesssim 2^{-Nd_3} |J|.$$
Using all these estimates, the term (\ref{cinu}) is seen to be smaller than

\begin{equation}
\sum_{d_1\geq 0} \sum_{K'\in \L_{d_1}} |K'|\sum_{d_1+d\leq d_3} 2^{d_3} 2^{d_3} 2^{-Nd_3} |J|=
\end{equation}

$$|J|\sum_{d_1\geq 0}2^{-(N-2)(d_1+d)} \sum_{K'\in \L_{d_1}} |K'|\lesssim$$

$$|I| |J| 2^{-(N-2)d} = 2^{(N-2)d}|R|,$$
by using (\ref{interval}), and this completes the proof.

\section{Counterexamples}

The next step in understanding this bi-parameter multi-linear framework is to consider more singular multipliers.
The most natural candidate is the double bilinear Hilbert transform, defined by

\begin{equation}\label{doubleBHT}
B_d(f,g)(x,y)=\int_{\R^2}f(x-t_1,y-t_2)g(x+t_1,y+t_2)\frac{dt_1}{t_1} \frac{dt_2}{t_2}=
\end{equation}

$$=\int_{\R^4} \sgn(\xi_1-\xi_2) \sgn(\eta_1-\eta_2)\widehat{f}(\xi_1, \eta_1)\widehat{g}(\xi_2, \eta_2)
e^{2\pi i(x,y)\cdot ((\xi_1,\eta_1) + (\xi_2,\eta_2))}\,d\xi d\eta.$$
It is the biparameter analogue of the bilinear Hilbert transform studied in \cite{lt} and given by

\begin{equation}\label{BHT}
B(f_1, f_2)(x)=\int_{\R}f_1(x-t) f_2(x+t) \frac{dt}{t}=
\end{equation}

$$\int_{\R^2} \sgn(\xi-\eta)\widehat{f}(\xi)\widehat{g}(\eta) e^{2\pi i x(\xi + \eta)}\, d\xi d\eta.$$
This time, the functions $f_1, f_2$ are defined on the real line. It is known (see \cite{lt}) that $B$ satisfies many $L^p$ estimates.

However, regarding $B_d$ we have the following theorem.

\begin{theorem}
The double bilinear Hilbert transform $B_d$ defined by (\ref{doubleBHT}), does not satisfy any $L^p$ estimates.
\end{theorem}

\begin{proof}
It is based on the following simple observation. Let $f(x,y)=g(x,y)=e^{ixy}$.
Since

$$(x-t_1)(y-t_2)+(x+t_1)(y+t_2)=2xy+2t_1 t_2$$
one can formally write

$$B(e^{ixy},e^{ixy})(x,y)=e^{2ixy}\int_{\R^2}e^{2it_1 t_2}\frac{dt_1}{t_1}\frac{dt_2}{t_2}=$$

$$4 e^{2ixy}\int_0^{\infty}\int_0^{\infty}\frac{\sin(t_1 t_2)}{t_1 t_2} dt_1 dt_2=$$

$$4e^{2ixy}\int_0^{\infty}(\int_0^{\infty}\frac{\sin(t_1 t_2)}{t_2} dt_2) \frac{dt_1}{t_1}=$$

$$4 e^{2ixy}\frac{\pi}{2}\int_0^{\infty} \frac{dt}{t}.$$

To obtain a quantitative version of this, we need the following lemma.

\begin{lemma}
There are two universal constants $C_1, C_2 >0$ such that 

\begin{equation}
\left|\int_0^N\int_0^N \frac{\sin(xy)}{xy} dx dy\right|\geq C_1 \log N
\end{equation}
as long as $N> C_2$.
\end{lemma}

\begin{proof}
Since $\int_0^{\infty}\frac{\sin t}{t} dt = \frac{\pi}{2}$, there is a constant
$C>0$ such that

\begin{equation}\label{bound}
\int_0^x\frac{\sin t}{t} dt \in [\frac{\pi}{4}, \frac{3\pi}{4}]
\end{equation}
whenever $x> C.$ Then,

$$\int_0^N\int_0^N \frac{\sin(xy)}{xy} dx dy=
\int_0^N(\int_0^N \frac{\sin(xy)}{y} dy) \frac{dx}{x}=
\int_0^N(\int_0^{Nx} \frac{\sin t}{t} dt) \frac{dx}{x}=$$

$$\int_0^{C/N}(\int_0^{Nx} \frac{\sin t}{t} dt) \frac{dx}{x}+
\int_{C/N}^N(\int_0^{Nx} \frac{\sin t}{t} dt) \frac{dx}{x}=$$

\begin{equation}\label{2terms}
\int_0^C(\int_0^{x} \frac{\sin t}{t} dt) \frac{dx}{x}+
\int_{C/N}^N(\int_0^{Nx} \frac{\sin t}{t} dt) \frac{dx}{x}.
\end{equation}
Since the function $x\rightarrow \frac{1}{x} \int_0^x \frac{\sin t}{t} dt$
is continuous on $[0,C]$ it follows that the first term in (\ref{2terms})
is actually $O(1)$. To estimate the second term in (\ref{2terms})
we observe that since $x> C/N$ it follows that $Nx > C$ and so, by using
(\ref{bound}) we can write

$$
\int_{C/N}^N(\int_0^{Nx} \frac{\sin t}{t} dt) \frac{dx}{x}\geq
\frac{\pi}{4}\int_{C/N}^N \frac{dx}{x} = \frac{\pi}{4}(2\log N - \log C),$$
and this ends the proof of the lemma, if $N$ is big enough.
\end{proof}

Now, coming back to the proof of the theorem, we define

$$f_N(x,y)= g_N(x,y)= e^{ixy} \chi_{[-N, N]}(x) \chi_{[-N, N]}(y)$$
and observe that

$$|B_d(f_N, g_N)(x,y)| \geq C
\left|\int_0^{N/10}\int_0^{N/10} \frac{\sin(zt)}{zt} dz dt\right|+
O(1)\geq C \log N + O(1)$$
as long as $x,y\in [-N/1000, N/1000]$.
This pointwise estimate precludes having $\|B_d(f_N, g_N)\|_r\leq C\|f_N\|_p \|g_N\|_q$ uniformly in $N$.

\end{proof}
At the end of this section, we would like to observe that, in the same manner,
 one can disprove the boundedness of the following operator
considered in \cite{cct2}.
Let $V$ be the trilinear operator $V$ defined by

\begin{equation}
V(f,g,h)(x)=\int_{\xi_1<\xi_2<\xi_3}\widehat{f}(\xi_1)
\widehat{g}(\xi_2)\widehat{h}(\xi_3) e^{2\pi i x(\xi_1-\xi_2+\xi_3)}\,
d\xi_1 d\xi_2 d\xi_3
\end{equation}

The following theorem holds (see \cite{cct2}).

\begin{theorem}
The trilinear operator $V$ constructed above does not map
$L^2\times L^2\times L^2\rightarrow L^{2/3,\infty}$.
\end{theorem}

\begin{proof}
First, by a simple change of variables one can reduce the study of $V$ to the study of $V_1$ defined by

\begin{equation}
V_1(f,g,h)(x)=\int_{\xi_1<-\xi_2<\xi_3}\widehat{f}(\xi_1)
\widehat{g}(\xi_2)\widehat{h}(\xi_3) e^{2\pi i x(\xi_1+\xi_2+\xi_3)}\,
d\xi_1 d\xi_2 d\xi_3.
\end{equation}
Also, we observe that the behaviour of $V_1$ is similar to the behaviour of 
$V_2$ defined by

\begin{equation}
V_2(f,g,h)(x)=\int_{\R^3}\sgn(\xi_1+\xi_2)\sgn(\xi_2+\xi_3)
\widehat{f}(\xi_1)\widehat{g}(\xi_2)\widehat{h}(\xi_3) 
e^{2\pi i x(\xi_1+\xi_2+\xi_3)}\,
d\xi_1 d\xi_2 d\xi_3,
\end{equation}
since the difference between $V_1$ and $V_2$ is a sum of simpler bounded operators.

But then, $V_2$ can be rewritten as

$$V_2(f,g,h)(x)=\int_{\R^2}f(x-t_1)g(x-t_1-t_2)h(x-t_2)\frac{dt_1}{t_1}
\frac{dt_2}{t_2}.$$
The counterexample is based on the following observation, similar to the one
 before.
Consider $f(x)=h(x)=e^{ix^2}$, $g(x)=e^{-ix^2}$. Because

$$(x-t_1)^2-(x-t_1-t_2)^2+(x-t_2)^2=x^2+2t_1t_2,$$
we can again formally write

$$V_2(e^{ix^2},e^{-ix^2}, e^{ix^2})(x)=e^{ix^2}\int_{\R^2}e^{2it_1t_2}\frac{dt_1}{t_1}\frac{dt_2}{t_2}=
4e^{ix^2}\frac{\pi}{4}\int_0^{\infty} \frac{dt}{t}.$$
To quantify this, we define
$f_N(x)= h_N(x)= e^{ix^2}\chi_{[-N, N]}(x)$ and 
$g_N(x)=e^{-ix^2}\chi_{[-N, N]}(x)$
and observe as before that

$$|V_2(f_N, g_N, h_N)(x)|\geq C
\left|\int_0^{N/10}\int_0^{N/10} \frac{\sin(xy)}{xy} dx dy\right|+
O(1)$$
if $x\in [-N/1000, N/1000]$ and this, as we have seen, contradicts 
the boundedness of the operator.

\end{proof}

\section{Further remarks}

First of all, we would like to remark that theorem (\ref{main}) has a straightforward generalization to the case
of $n$-linear operators, for $n\geq 1$. 

Let $m\in L^{\infty}(\R^{2n})$ be a symbol satisfying the bi-parameter Marcinkiewicz-H\"{o}rmander-Mihlin condition

\begin{equation}\label{bibi}
|\partial_{\xi}^{\alpha}\partial_{\eta}^{\beta} m(\xi, \eta)|\lesssim \frac{1}{|\xi|^{|\alpha|}}
\frac{1}{|\eta|^{|\beta|}},
\end{equation}
for many multiindices $\alpha$ and $\beta$.
Then, for $f_1,...,f_n$ Schwartz functions in $\R^2$, define the operator $T_m$ by

\begin{equation}\label{ultimul}
T_m(f_1,...,f_n)(x):= \int_{\R^{2n}} m(\xi,\eta) \widehat{f_1}(\xi_1,\eta_1)...\widehat{f_n}(\xi_n,\eta_n)
e^{2\pi i x\cdot ((\xi_1,\eta_1)+...+(\xi_n,\eta_n))}\,d\xi d\eta.
\end{equation}
We thus record

\begin{theorem}
The bi-parameter $n$-linear operator $T_m$ maps $L^{p_1}\times ... \times L^{p_n}\rightarrow L^p$
as long as $1<p_1,...,p_n\leq\infty$, $1/p_1 + ... + 1/p_n = 1/p$ and $0<p<\infty$.
\end{theorem}

Here, when such an $n+1$-tuple $(p_1,...,p_n,p)$ has the property that $0<p<1$ and $p_j=\infty$ for some $1\leq j\leq n$
then, for some technical reasons (see \cite{cct1}), by $L^{\infty}$ one actually means $L^{\infty}_c$ the space of bounded measurable functions with compact support.

On the other hand, one can ask what is happening if one is interested in more singular multipliers. Suppose $\Gamma_1$ and $\Gamma_2$ are subspaces in $\R^n$ and one 
considers operators $T_m$ defined by (\ref{ultimul})
where $m$ satisfies

\begin{equation}\label{gamas}
|\partial_{\xi}^{\alpha}\partial_{\eta}^{\beta} m(\xi, \eta)|\lesssim \frac{1}{\dist(\xi,\Gamma_1)^{|\alpha|}}
\frac{1}{|\dist(\eta,\Gamma_2)^{|\beta|}}.
\end{equation}
Our theorem says that if $\dim(\Gamma_1)=\dim(\Gamma_2)=0$ then we have many $L^p$ estimates available. 
On the other hand,
the previous counterexamples show that  when $\dim(\Gamma_1)=\dim(\Gamma_2)=1$ then we do not have any $L^p$ estimates.
But it is of course natural to ask 

\begin{question}
Let $\dim(\Gamma_1)=0$ and $\dim(\Gamma_2)=1$ with $\Gamma_2$ non-degenerate in the sense of \cite{cct1}. If $m$ is a multiplier satisfying (\ref{gamas}) does the corresponding $T_m$ satisfy any $L^p$ estimates ?
\end{question}

\section{Appendix: differentiating paraproducts }

In this section we describe how the Kato-Ponce inequality (\ref{katoponce}) can be reduced to Coifman-Meyer theorem
and also how the more general inequality (\ref{2katoponce}) can be reduced to our theorem \ref{main}.

The argument is standard and is based on some "calculus with paraproducts". 

In what follows, we will define generic classes of paraproducts. First we consider the sets $\Phi$ and
$\Psi$ given by

$$\Phi:= \{ \phi\in\S(\R) : \supp\widehat{\phi}\subseteq [-1,1] \},$$

$$\Psi:= \{ \psi\in\S(\R) : \supp\widehat{\psi}\subseteq [1,2] \}.$$
The intervals $[-1,1]$ and $[1,2]$ are not important. What is important, is the fact that $\Phi$ consists of
Schwartz functions whose Fourier support is compact and contains the origin and $\Psi$ consists of Schwartz functions whose Fourier support is compact and does not contain the origin. Then, for various $\phi\in\Phi$ and
$\psi,\psi'\psi''\in\Psi,$ we define the paraproducts $\Pi_j$ $j=0,1,2,3$ as follows

\begin{equation}\label{pi0}
\Pi_0(f,g)(x):=\int_{\R}\left((f\ast D^1_{2^k}\psi)(g\ast D^1_{2^k}\psi')\right)\ast D^1_{2^k}\psi''(x)\, dk,
\end{equation}

\begin{equation}\label{pi1}
\Pi_1(f,g)(x):=\int_{\R}\left((f\ast D^1_{2^k}\phi)(g\ast D^1_{2^k}\psi)\right)\ast D^1_{2^k}\psi'(x)\, dk,
\end{equation}

\begin{equation}\label{pi2}
\Pi_2(f,g)(x):=\int_{\R}\left((f\ast D^1_{2^k}\psi)(g\ast D^1_{2^k}\phi)\right)\ast D^1_{2^k}\psi'(x)\, dk,
\end{equation}

\begin{equation}\label{pi3}
\Pi_3(f,g)(x):=\int_{\R}\left((f\ast D^1_{2^k}\psi)(g\ast D^1_{2^k}\psi')\right)\ast D^1_{2^k}\phi(x)\,dk.
\end{equation}
All these paraproducts are bilinear operators for which the Coifman-Meyer theorem applies. For instance, one can rewrite
$\Pi_0(f,g)$ as

$$\Pi_0(f,g)(x)= \int_{\R^2} m(\xi_1,\xi_2) \widehat{f}(\xi_1)\widehat{g}(\xi_2) e^{2\pi i x(\xi_1+\xi_2)}\, d\xi_1 d\xi_2,$$
where the symbol $m(\xi_1 ,\xi_2)$ is given by

$$m(\xi_1,\xi_2)= \int_{\R}
(D^{\infty}_{2^{-k}}\widehat{\psi})(\xi_1)
(D^{\infty}_{2^{-k}}\widehat{\psi'})(\xi_2)
(D^{\infty}_{2^{-k}}\widehat{\psi''})(-\xi_1-\xi_2)\, dk,$$
and satisfies the Marcinkiewicz-H\"{o}rmander-Mihlin condition.

The reduction relies on the follwing simple observation

\begin{proposition}\label{calculus}
Let $\alpha>0$. Then, for every paraproduct $\Pi_1$ there exists a paraproduct $\Pi'_1$ so that

\begin{equation}
\D^{\alpha}\Pi_1(f,g) = \Pi'_1(f,\D^{\alpha}g),
\end{equation}
for every $f,g$ Schwartz functions on $\R$.
\end{proposition}
\begin{proof}
It is based on the following equalities

$$\D^{\alpha}\Pi_1(f,g)=$$

$$\int_{\R}\left((f\ast D^1_{2^k}\phi)(g\ast D^1_{2^k}\psi)\right)\ast \D^{\alpha}(D^1_{2^k}\psi')\, dk$$

$$\int_{\R}\left((f\ast D^1_{2^k}\phi)(g\ast D^1_{2^k}\psi)\right)\ast 2^{-k\alpha}D^1_{2^k}(\D^{\alpha}\psi')\, dk=$$

$$\int_{\R}\left((f\ast D^1_{2^k}\phi)(g\ast 2^{-k\alpha}D^1_{2^k}\psi)\right)\ast D^1_{2^k}(\D^{\alpha}\psi')\, dk=$$

$$\int_{\R}\left((f\ast D^1_{2^k}\phi)(g\ast \D^{\alpha}(D^1_{2^k}(\D^{-\alpha}\psi)))\right)\ast D^1_{2^k}(\D^{\alpha}\psi')\, dk=$$

$$\int_{\R}\left((f\ast D^1_{2^k}\phi)(\D^{\alpha}g\ast D^1_{2^k}(\D^{-\alpha}\psi))\right)\ast D^1_{2^k}(\D^{\alpha}\psi')\, dk:=$$

$$\Pi'_1(f,\D^{\alpha}g),$$
where $\D^{-\alpha}\psi$ is the Schwartz function whose Fourier transform is given by 
$\widehat{\D^{-\alpha}\psi}(\xi)=|\xi|^{-\alpha}\widehat{\psi}(\xi)$, which is well defined since $\psi\in\Psi$.
\end{proof}

Clearly, one has similar identities for all the other types of paraproducts $\Pi_j$ $j\neq 1$. However, one has to be particularly careful about the case
of $\Pi_3$ since there, the corresponding functions $|\xi|^{\alpha}\widehat{\phi}(\xi)$ are no longer smooth, and as a consequence, their inverse Fourier transforms
have only limited decay of type $1/(1+|x|)^{1+\alpha}$.
To prove the Kato-Ponce inequality (\ref{katoponce}), one just has to realize that every product of two functions $f$ and $g$ on $\R$, 
can be written as a sum of such paraproducts

$$fg=\sum_{j=0}^3 \Pi_j(f,g)$$
and then, after using the above proposition \ref{calculus}, to apply the Coifman-Meyer theorem. In fact, the argument of this paper can be naturally strenghten to prove more,
namely that (\ref{katoponce}) holds as long as $1<p,q\leq \infty$, $1/r=1/p+1/q$ and $1/(1+\alpha) < r < \infty$. See the second volume of \cite{ms} for details.
The constraint $1/(1+\alpha) < r < \infty$ is clearly related to the limited decay mentioned above.

A similar treatment is available in the bi-parameter case too. Here, one has to handle bi-parameter paraproducts 
$\Pi_{i,j}$ for $i,j=0,1,2,3$ formally defined by $\Pi_{i,j}:=\Pi_i\otimes\Pi_j$.

One first observes the following extension of proposition \ref{calculus}

\begin{proposition}\label{2calculus}
Let $\alpha$, $\beta>0$. Then, for every paraproduct $\Pi_{1,2}$ there exists a paraproduct $\Pi'_{1,2}$ so that
\begin{equation}
\D_1^{\alpha}\D_2^{\beta}\Pi_{1,2}(f,g) = \Pi'_{1,2}(\D_2^{\beta}f,\D_1^{\alpha}g),
\end{equation}
for every $f,g$ Schwartz functions on $\R^2$.
\end{proposition}

There are also similar equalities for the remaining paraproducts $\Pi_{i,j}$ when $(i,j)\neq (1,2)$. 
Since every product of two functions $f$ and $g$ on $\R^2$ can be written as

$$fg=\sum_{i,j=0}^3 \Pi_{i,j}(f,g),$$
everything follows from theorem \ref{main}. In fact, as before, the techniques of the present paper can be strenghten and one can similarly prove that an
even more general inequality holds, namely

\begin{equation}
\|\D_1^{\alpha}\D_2^{\beta}(fg)\|_r \lesssim  \|\D_1^{\alpha}\D_2^{\beta}f\|_{p_1} \|g\|_{q_1} + \|f\|_{p_2} \|\D_1^{\alpha}\D_2^{\beta}g\|_{q_2} + 
 \|\D_1^{\alpha}f\|_{p_3} \|\D_2^{\beta}g\|_{q_3} + \|\D_1^{\alpha}g\|_{p_4} \|\D_2^{\beta}f\|_{q_4}
\end{equation}
whenever $1<p_j, q_j\leq\infty$, $1/p_j + 1/q_j =1/r$ for $j=1,2,3,4$ and $\max (1/(1+\alpha), 1/(1+\beta)) <r<\infty$. See again the second volume of the book
\cite{ms} for details.

Finally, we would like to mention that when the first version of this article has been released, we have not been particularly careful about the precise conditions
under which these more general forms of the Kato-Ponce inequalities hold (both in the one-parameter and multi-parameter case). 
We would like to thank
Loukas Grafakos and Seungly Oh for pointing this oversight to us. In fact, their recent paper \cite{go} obtains independently, similar estimates.

\end{document}